\vsize=9.0in\voffset=1cm
\looseness=2


\message{fonts,}

\font\tenrm=cmr10
\font\ninerm=cmr9
\font\eightrm=cmr8
\font\teni=cmmi10
\font\ninei=cmmi9
\font\eighti=cmmi8
\font\ninesy=cmsy9
\font\tensy=cmsy10
\font\eightsy=cmsy8
\font\tenbf=cmbx10
\font\ninebf=cmbx9
\font\tentt=cmtt10
\font\ninett=cmtt9

\font\ninesl=cmsl9
\font\eightsl=cmsl8

\font\nineit=cmti9
\font\eightit=cmti8

\skewchar\ninei='177 \skewchar\eighti='177
\skewchar\ninesy='60 \skewchar\eightsy='60

\def\eightpoint{\def\rm{\fam0\eightrm} 
\normalbaselineskip=9pt
\normallineskiplimit=-1pt
\normallineskip=0pt

\textfont0=\eightrm \scriptfont0=\sevenrm \scriptscriptfont0=\fiverm
\textfont1=\ninei \scriptfont1=\seveni \scriptscriptfont1=\fivei
\textfont2=\ninesy \scriptfont2=\sevensy \scriptscriptfont2=\fivesy
\textfont3=\tenex \scriptfont3=\tenex \scriptscriptfont3=\tenex
\textfont\itfam=\eightit  \def\it{\fam\itfam\eightit} 
\textfont\slfam=\eightsl \def\sl{\fam\slfam\eightsl} 

\setbox\strutbox=\hbox{\vrule height6pt depth2pt width0pt}%
\normalbaselines \rm}

\def\ninepoint{\def\rm{\fam0\ninerm} 
\textfont0=\ninerm \scriptfont0=\sevenrm \scriptscriptfont0=\fiverm
\textfont1=\ninei \scriptfont1=\seveni \scriptscriptfont1=\fivei
\textfont2=\ninesy \scriptfont2=\sevensy \scriptscriptfont2=\fivesy
\textfont3=\tenex \scriptfont3=\tenex \scriptscriptfont3=\tenex
\textfont\itfam=\nineit  \def\it{\fam\itfam\nineit} 
\textfont\slfam=\ninesl \def\sl{\fam\slfam\ninesl} 
\textfont\bffam=\ninebf \scriptfont\bffam=\sevenbf
\scriptscriptfont\bffam=\fivebf \def\bf{\fam\bffam\ninebf} 
\textfont\ttfam=\ninett \def\tt{\fam\ttfam\ninett} 

\normalbaselineskip=11pt
\setbox\strutbox=\hbox{\vrule height8pt depth3pt width0pt}%
\let \smc=\sevenrm \let\big=\ninebig \normalbaselines
\parindent=1em
\rm}

\def\tenpoint{\def\rm{\fam0\tenrm} 
\textfont0=\tenrm \scriptfont0=\ninerm \scriptscriptfont0=\fiverm
\textfont1=\teni \scriptfont1=\seveni \scriptscriptfont1=\fivei
\textfont2=\tensy \scriptfont2=\sevensy \scriptscriptfont2=\fivesy
\textfont3=\tenex \scriptfont3=\tenex \scriptscriptfont3=\tenex
\textfont\itfam=\nineit  \def\it{\fam\itfam\nineit} 
\textfont\slfam=\ninesl \def\sl{\fam\slfam\ninesl} 
\textfont\bffam=\ninebf \scriptfont\bffam=\sevenbf
\scriptscriptfont\bffam=\fivebf \def\bf{\fam\bffam\tenbf} 
\textfont\ttfam=\tentt \def\tt{\fam\ttfam\tentt} 

\normalbaselineskip=11pt
\setbox\strutbox=\hbox{\vrule height8pt depth3pt width0pt}%
\let \smc=\sevenrm \let\big=\ninebig \normalbaselines
\parindent=1em
\rm}

\message{fin format jgr}

\vskip 4 mm
\magnification=1200
\font\Bbb=msbm10
\def\R{\hbox{\Bbb R}}

\def\S{\hbox{\Bbb S}}

\def\b{\backslash}

\vskip 4 mm

\centerline{\bf Weighted Radon transforms for which}
\centerline{\bf the Chang approximate inversion formula is precise}

\centerline{\bf R.G. Novikov}
\vskip 4 mm
\noindent
{\ninerm CNRS (UMR 7641), Centre de Math\'ematques Appliqu\'ees, Ecole
Polytechnique,}

\noindent
{\ninerm 91128 Palaiseau, France and}

\noindent
{\ninerm IIEPT RAS - MITPAN, Profsoyuznaya str., 84/32,  Moscow 117997,
Russia}

\noindent
{\ninerm e-mail: novikov@cmap.polytechnique.fr}

\vskip 4 mm
{\bf Abstract.}\
We describe all weighted Radon transforms on the plane for which
the Chang approximate inversion formula is precise. Some subsequent
results, including  the Cormack type inversion for these transforms, are also
given.

\vskip 2 mm
{\bf 1.Introduction}

We consider the weighted ray transformation $P_W$ defined by the formula
$$\eqalign{
&P_Wf(s,\theta)=\int\limits_{\R} W(s\theta^{\perp}+t\theta,\theta)
f(s\theta^{\perp}+t\theta)dt,\cr
&s\in\R,\ \theta=(\theta_1,\theta_2)\in\S^1,\
\theta^{\perp}=(-\theta_2,\theta_1),\cr}\eqno(1)$$
where $W=W(x,\theta)$ is the weight, $f=f(x)$ is a test function, $x\in\R$,
$\theta\in\S^1$.
Up to change of variables, $P_W$ is known also as weighted Radon transform
on the plane.

We recall that in definition (1) the product $\R\times\S^1$ is interpreted as
the set of all oriented straight lines in $\R^2$. If
$\gamma=(s,\theta)\in\R\times\S^1$, then
$\gamma=\{x\in\R^2:\ x=s\theta^{\perp}+t\theta,\ t\in\R\}$ (modulo
orientation) and $\theta$ gives the orientation of $\gamma$.

We assume that
$$\eqalign{
&W\ \ {\rm is\ complex-valued},\cr
&W\in C(\R^2\times\S^1)\cap L^{\infty}(\R^2\cap\S^1),\cr
&w_0(x)\buildrel \rm def \over =
{1\over 2\pi}\int\limits_{\S^1}W(x,\theta)d\theta\ne 0,\ \ x\in\R^2,\cr}
\eqno(2)$$
where $d\theta$ is the standard element of arc length on $\S^1$.

If $W\equiv 1$, then $P_W$ is known as the classical ray (or Radon) transform
on the plane. If
$$\eqalign{
&W(x,\theta)=\exp(-D a(x,\theta)),\cr
&D a(x,\theta)=\int\limits_0^{+\infty}a(x+t\theta)dt,\cr}\eqno(3)$$
where $a$ is a complex-valued sufficiently regular function on $\R^2$ with
sufficient decay at infinity, then $P_W$ is known as the attenuated ray
(or Radon) transform.

The classical Radon transform arises, in particular, in the X-ray transmission
 tomography. The attenuated Radon transform (at least, with $a\ge 0$) arises,
in particular, in the single photon emission computed tomography (SPECT).
Some other weights $W$ also arise in applications. For more information in
this connection see, for example, [Na], [K].

Precise and simultaneously explicite inversion formulas for the classical
and attenuated Radon transforms were given for the first time in [R] and
[No], respectively. For some other weights $W$ precise and
simultaneously explicite inversion formulas  were given in [BS], [G].

On the other hand, the following Chang approximate inversion formula for
$P_W$, where $W$ is given by (3) with $a\ge 0$, is used for a long time,
see [Ch], [M], [K]:
$$\eqalign{
&f_{appr}(x)={1\over 4\pi w_0(x)}\int\limits_{\S^1}
h^{\prime}(x\theta^{\perp},\theta)d\theta,\ \
h^{\prime}(s,\theta)={d\over ds}h(s,\theta),\cr
&h(s,\theta)={1\over \pi}p.v.\int\limits_{\R}
{P_Wf(t,\theta)\over {s-t}}dt,\ \ s\in\R,\ \theta\in\S^1,\ x\in\R^2,\cr}
\eqno(4)$$
where $w_0$ is defined in (2). It is known that (4) is efficient as the
first approximation in SPECT reconstructions and, in particular, is
sufficiently stable to the strong Poisson noise of SPECT data.
The  results of the present note consist of the following:
\item{(1)} In Theorem, under assumptions (2), we describe all weights $W$
for which the Chang approximate inversion formula (4) is precise, that is
$f_{appr}\equiv f$ on $\R^2$;
\item{(2)} For $P_W$ with $W$ of Theorem  we give also the Cormack type
inversion (see Remark A) and inversion from limited angle data
(see Remark B).

These results are presented in detail in the next section. In addition, we
give also an explanation of  efficiency of the  Chang
formula (4) as the first approximation in SPECT reconstructions
(on the level of integral geometry).

\vskip 2 mm
{\bf 2. Results}

Let
$$\eqalign{
&C_0(\R^2)\ \ {\rm denote\ the\ space\ of\ continuous}\cr
&{\rm compactly\ supported\ functions\ on}\ \ \R^2.\cr}\eqno(5)$$
Let
$$\eqalign{
&L^{\infty,\sigma}(\R^2)=\{f:\ M^{\sigma}f\in L^{\infty}(\R^2)\},\cr
&M^{\sigma}f(x)=(1+|x|)^{\sigma}f(x),\ \ x\in\R^2,\ \ \sigma\ge 0.\cr}\eqno(6)
$$

\vskip 2 mm
{\bf Theorem.}
{\sl
Let assumptions (2) hold and let $f_{appr}(x)$ be given by (4). Then
$$f_{appr}=f\ \ { (in\ the\ sense\ of\ distributios)\ on}\ \ \R^2\ \
{ for\ all}\ \ f\in C_0(\R^2) \eqno(7)$$
if and only if
$$W(x,\theta)-w_0(x)\equiv w_0(x)-W(x,-\theta),\ \ x\in\R^2,\ \ \theta\in\S^1.
\eqno(8)$$
(This result remains valid with $C_0(\R^2)$ replaced by
$L^{\infty,\sigma}(\R^2)$ for $\sigma>1$.)
}

Theorem 1 is based on the following facts:
\item{$\bullet$} Formula (4) coincides with the classical Radon inversion formula
if $W\equiv 1$ and, as a corollary, is precise if $W\equiv w_0$.
\item{$\bullet$} Formula (4) is equivalent to the symmetrized formula
$$\eqalign{
&f_{appr}(x)={1\over 4\pi w_0(x)}\int\limits_{\S^1}
g^{\prime}(x\theta^{\perp},\theta)d\theta,\ x\in\R^2,\cr
&g(s,\theta)={1\over 2\pi}p.v.\int\limits_{\R}
{{P_Wf(t,\theta)+P_Wf(-t,-\theta)}\over {s-t}}dt,\ \
(s,\theta)\in\R\times\S^1.\cr}\eqno(9)$$
\item{$\bullet$} The following formula holds:
$$\eqalign{
&{1\over 2}(P_Wf(s,\theta)+P_Wf(-s,-\theta))=P_{W_{sym}}f(s,\theta),
(s,\theta)\in\R\times\S^1,\cr
&W_{sym}(x,\theta)={1\over 2}(W(x,\theta)+W(x,-\theta)),\ \ x\in\R^2,\
\theta\in\S^1.\cr}\eqno(10)$$
\item{$\bullet$} If
$$\eqalign{
&q\in C(\R\times\S^1),\  supp\,q\ \ {\rm is\ compact},\ \
q(s,\theta)=q(-s,-\theta),\cr
&g(s,\theta)={1\over \pi}p.v.\int\limits_{\R}{q(t,\theta)\over {s-t}}dt,\ \
(s,\theta)\in\R\times\S^1,\cr
&\int\limits_{\S^1}g^{\prime}(x\theta^{\perp},\theta)d\theta\equiv 0\ \
{\rm as\ a\ distribution\ of}\ \ x\in\R^2,\cr}\eqno(11)$$
then $q\equiv 0$ on $\R\times\S^1$.

The statement that, under the assumptions of Theorem, property (8) implies
(7) can be also deduced from considerations developed in [K].

Using that $W_{sym}\equiv w_0$ under condition (8), we obtain also the
following

\vskip 2 mm
{\bf Remarks.}\
Let conditions (2), (8) be fulfiled. Let $f\in C_0(\R^2)$. Then:
\item{(A)} $P_Wf$ on $\Omega(D)$ uniquely determines $f$ (or more precisely
$w_0f$) on $\R^2\b D$ via (10) and the Cornack inversion from
$P_{W_{sym}}f$ on $\Omega(D)$, where $D$ is a compact in $\R^2$,
$\Omega(D)$ denotes the set of  all straight lines in $\R^2$ which do not
intersect $D$;
\item{(B)}  $P_Wf$ on $\R\times (S_+\cup S_-)$ uniquely determines $f$ on
$\R^2$ via (10) and standard inversion from the limited angle data
$P_{W_{sym}}f$ on $\R\times S_+$, where $S_+$ is an arbitrary nonempty open
connected subset of $\S^1$, $S_-=-S_+$.

For the case when $W$ is given by  (3) under the additional conditions that
$a\ge 0$ and $supp\,a\subset D$, where $D$ is some known bounded domain
which is not too big, and for $f\in C(\R^2)$, $f\ge 0$, $supp\,f\subset D$,
the transform $P_Wf$ is relatively well approximated by $P_{W_{appr}}f$,
where $W_{appr}(x,\theta)=w_0(x)+(1/2)(W(x,\theta)-W(x,-\theta))$.
In addition,  this $W_{appr}$ already satisfies (8). This explains the
efficiency of (4) as the first aproximation in SPECT reconstructions
(on the level of integral geometry).

\vfill\eject
\vskip 4 mm
{\bf References}
\item{[BS]} J.Boman  and J.O.Str\"omberg,  Novikov's inversion formula for
the attenuated Radon transform - a new approach, {\it J.Geom.Anal.} {\bf 14}
(2004), 185-198
\item{[Ch]} L.T.Chang, A method for attenuation correction in
radionuclide computed tomography,
{\it IEEE Trans. Nucl. Sci.} NS-25 (1978), 638-643
\item{[ G]} S.Gindikin, A remark on the weighted Radon transform on the
plane, {\it Inverse Problems and Imaging} {\bf 4} (2010), 649-653
\item{[ K]} L.A.Kunyansky, Generalized and attenuated Radon transforms:
restorative approach to the numerical inversion, {\it Inverse Problems}
{\bf 8} (1992), 809-819
\item{[ M]} K.Murase, H.Itoh, H.Mogami, M.Ishine, M.Kawamura, A.Lio  and
K.Hamamoto,  A comparative study of attenuation correction algorithms in
single photon emission computed tomography (SPECT), {\it Eur. J. Nucl. Med.}
{\bf 13} (1987), 55-62
\item{[Na]} F.Natterer, {\it The Mathematics of Computerized Tomography}
(Stuttgart: Teubner), 1986
\item{[No]} R.G.Novikov, An inversion formula for the attenuated X-ray
transformation, {\it Ark. Mat.} {\bf 40} (2002), 145-167
\item{[ R]} J.Radon, Uber die Bestimmung von Funktionen durch ihre
Integralwerte langs

\item{    } gewisser Mannigfaltigkeiten, {\it Ber. Verh. Sachs. Akad.
Wiss. Leipzig, Math-Nat.}, {\bf K}1 {\bf 69} (1917), 262-267

\end